\documentclass[a4paper,12pt]{amsart}
\usepackage{amssymb}

\usepackage[colorlinks=true, linktocpage=false]{hyperref}
\numberwithin{equation}{section}
\newtheorem{theorem}{Theorem}[section]
\newtheorem{prop}[theorem]{Proposition}
\newtheorem{lem}[theorem]{Lemma}

\theoremstyle{remark}
\newtheorem{rem}[theorem]{Remark}
\newcommand{\R}{\mathbb{R}}

\newcommand{\N}{\mathbb{N}}

\author[J.~Benameur]{Jamel Benameur}
\author[M.~Ltifi]{Maroua Ltifi}
\address{university of Gab\`es, Faculty of Science of Gab\`es, Department of Mathematics; Tunisia}
\email{\sl jamelbenameur@gmail.com}
\email{\sl widaltifi@gmail.com}

\title[Strong solution of 3D-NSE  with exponential damping]
{Strong solution of 3D-NSE  with exponential damping}

\begin{document}
\begin{abstract}
In this paper we prove the existence and uniqueness of strong
solution of the incompressible Navier-Stokes equations with
damping $\alpha (e^{\beta|u|^2}-1)u$.
\end{abstract}



\maketitle
\tableofcontents


\section{\bf Introduction}
    The classical Navier-Stokes equation is an area wish has received some attention during the last period, where as
    the study of this equation has become classical see \cite{HB},\cite{JYC}. Despite, the solution of 3D Navier-Stokes
    equations is still a big open problem although. In 2008, the modified Navier-Stokes equations with
    damping $\alpha|u|^{\beta-1}u$, was studied by Cai and Jiu \cite{XJ}, they proved the global existence of a weak solution if $u^0$ in $L^2(\R^3)$ and they proved the global existence and uniqueness of a strong solution if the initial condition $u^0$ is in $H^1(\R^3)\cap L^{\beta+1}(\R^3)$, with $\beta\geq7/2$. To construct a global solution Cai and Jiu used the Galerkin approximations. There is a large of literature
    dealing the classical Navier-Stokes equations in different spaces. Recently, Benameur in \cite{JB} has considered a new model of the
    Navier-stokes equation called  Navier-Stokes equations with exponential damping, where he proved the global existence of weak solution. In this paper we study
    the global existence of strong solution to the incompressible Navier-Stokes equations
    with exponential damping in three spatial dimensions
$$(NS)
  \begin{cases}
     \partial_t u
 -\nu_h\Delta_h u-\nu_3\partial_3^2u+ u.\nabla u  +\alpha (e^{\beta |u|^2}-1)u =\;\;-\nabla p\hbox{ in } \mathbb R^+\times \mathbb R^3\\
     {\rm div}\, u = 0 \hbox{ in } \mathbb R^+\times \mathbb R^3\\
    u(0,x) =u^0(x) \;\;\hbox{ in }\mathbb R^3,
  \end{cases}
$$
where $\nu_h>0$ and $\nu_3\geq0$ are respectively the horizontal and vertical viscosity of fluid,
$u=u(t,x)=(u_1,u_2,u_3)$ and $p=p(t,x)$ denote respectively the unknown velocity and the unknown
pressure of the fluid at the point $(t,x)\in \mathbb R^+\times \mathbb R^3$, and $\alpha,\beta>0$.
The terms $(u.\nabla u):=u_1\partial_1 u+u_2\partial_2 u+u_3\partial_3u$, while
$u^0=(u_1^o(x),u_2^o(x),u_3^o(x))$ is an initial given velocity. If $u^0$ is quite regular,
the divergence free condition determines the pressure $p$. First, we study the isotropic case $\nu_h=\nu_3=1$:
$$(NS_1)
  \begin{cases}
     \partial_t u
 -\Delta u+ u.\nabla u  +\alpha (e^{\beta |u|^2}-1)u =\;\;-\nabla p\hbox{ in } \mathbb R^+\times \mathbb R^3\\
     {\rm div}\, u = 0 \hbox{ in } \mathbb R^+\times \mathbb R^3\\
    u(0,x) =u^0(x) \;\;\hbox{ in }\mathbb R^3,
  \end{cases}
$$
and the associated space:
$$\begin{array}{lcl}\mathcal H_{\beta}&=&\{f:\mathbb R^+\times\mathbb R^3\rightarrow\mathbb
R^3\,measurable;\\
&&(e^{\beta |f|^2}-1)|f|^2,(e^{\beta |f|^2}-1)|\nabla
f|^2,\;e^{\beta |f|^2}|\nabla |f|^2|^2\in
L^1(\R^+\times\R^3)\}.\end{array}$$ Our first result is the
following.
\begin{theorem}\label{th1} Let $u^0\in H^1(\mathbb R^3)$ be a divergence free vector fields, then there is a unique global solution of $(NS_1)$:
    $u\in L^\infty(\R^+,H^1(\mathbb R^3)\cap C(\R^+,H^{-2}(\R^3))\cap L^2(\R^+,\dot H^2(\mathbb R^3))\cap \mathcal H_{\beta}$. Moreover, for all $t\geq0$
\begin{equation}\label{eqth1}\|u(t)\|_{L^2}^2+2\int_0^t\|\nabla u\|_{L^2}^2+2\alpha\int_0^t\|(e^{\beta |u|^2}-1)|u|^2\|_{L^1}\leq \|u^0\|_{L^2}^2,\end{equation}
\begin{equation}\label{eqth2}\|\nabla u(t)\|_{L^2}^2+\int_0^t\|\Delta u\|_{L^2}^2+\alpha \beta\int_{0}^{t}\|e^{\beta|u|^{2}}|\nabla(|u|^{2})|^{2}\|_{L^{1}}+\alpha\int_0^t\|(e^{\beta |u|^2}-1)|\nabla u|^2\|_{L^1}\leq \|\nabla u^0\|_{L^2}^2e^{\frac{t}{\alpha\beta^2}},\end{equation}
\begin{equation}\label{eqth3}\|\nabla u(t)\|_{L^2}^2+\int_0^t\|\Delta u\|_{L^2}^2+\alpha \beta\int_{0}^{t}\|e^{\beta|u|^{2}}(\nabla(|u|^{2})|^{2}\|_{L^{1}}\\
+\alpha\int_0^t\|(e^{\beta |u|^2}-1)|\nabla u|^2\|_{L^1}\leq M_{\alpha,\beta}(u^0),
\end{equation}
where $M_{\alpha,\beta}(u^0)=\|\nabla u^0\|_{L^2}^2+\frac{\|u^0\|_{L^2}^2}{\alpha\beta^2}.$
\end{theorem}
\begin{rem}\label{rem1}\begin{enumerate}
\item The fact $(e^{\beta |u|^2}-1)|u|^2\in L^1(\R^+,L^1(\R^3))$
implies $u\in \cap _{4\leq p<\infty}L^p(\R^+,L^p(\R^3))$. Indeed:
we have
$$\int_0^\infty\|(e^{\beta|u(t)|^2}-1)|u(t)|^2\|_{L^1}dt=\sum_{k=4}^{\infty}\frac{\beta^k}{k!}\int_0^\infty\|u(t)\|_{L^{2k+2}}^{2k+2}dt.$$
\item By interpolation between $H^{-2}(\R^3)$ and $H^1(\R^3)$, we
obtain: For all $s<1$, we have $u\in C(\mathbb R^+,H^{s}(\R^3))$.
\item $u\in C_r(\mathbb R^+,H^1(\R^3))$. Indeed: By equations
(\ref{eqth1})-(\ref{eqth2}) we get
$\limsup_{t\rightarrow0^+}\|u(t)\|_{H^1}\leq \|u^0\|_{H^1}$.
Applying Proposition \ref{prop1}, we get the continuity of $u$ at
0. For $t_0>0$, consider the following system
$$
  (S(t_0))\begin{cases}
     \partial_t v
 -\Delta v+ v.\nabla v  +\alpha (e^{\beta |v|^2}-1)v =\;\;-\nabla q\hbox{ in } \mathbb R^+\times \mathbb R^3\\
     {\rm div}\, v = 0 \hbox{ in } \mathbb R^+\times \mathbb R^3\\
    v(0,x) =u(t_0,x) \;\;\hbox{ in }\mathbb R^3.
  \end{cases}
$$
By the uniqueness given by Theorem \ref{th1}, we obtain
$v(t,x)=u(t_0+t,x)$ is the global solution of $(S(t_0))$. Then $v$
is right continuous at 0, which implies the right continuity of
$u$ at $t_0$.
\item The continuity set of $u$ solution of $(NS_1)$ in $H^1(\R^3)$: Put the following subset of $\R^+$
$$A=\{t\in\R^+;\;u ~~~\mbox{discontinous at} ~~~t~~~\mbox{in}\,H^1(\R^3)\}.$$
$A$ is at most countable. Indeed: Let $f(t)=e^{\frac{-t}{\alpha\beta^{2}}}\|\nabla u(t)\|^{2}_{L^2}$, $g(t)=\|\nabla u(t)\|^{2}_{L^2}$ and
 $$B=\{t\in\R^+;\;f ~~~\mbox{discontinous at} ~~~t\}=\{t\in\R^+;\;g ~~~\mbox{discontinous at} ~~~t\}.$$
Let $0\leq t_{1}< t_{2}$. Combinig the uniqueness of strong solution of $(NS_1)$ and inequality (\ref{eqth2}), we get
$$ \|\nabla u(t_{2})\|^{2}_{L^2}\leq\|\nabla u(t_{1})\|^{2}_{L^2} e^{\frac{t_{2}-t_{1}}{\alpha\beta^{2}}}$$
and
  $$\|\nabla u(t_{2})\|^{2}_{L^2}e^{\frac{-t_{2}}{\alpha\beta^{2}}}\leq \|\nabla u(t_{1})\|^{2}_{L^2} e^{\frac{-t_{1}}{\alpha\beta^{2}}}.$$
Thus, $f$ is a decreasing function. According  Lemma \ref{lem5}, $B$ is at most countable. By using Proposition \ref{prop1} and Remark \ref{remprop1}, we obtain $B=A$ and the desired result is proved.
\end{enumerate}
\end{rem}

Secondly, we study the anisotropic Navier-Stokes case $\nu_h=1,\,\nu_3=0$ with the same damping:
$$(NS_2)
  \begin{cases}
     \partial_t u
 -\Delta_h u+ u.\nabla u  +\alpha (e^{\beta |u|^2}-1)u =\;\;-\nabla p\hbox{ in } \mathbb R^+\times \mathbb R^3\\
     {\rm div}\, u = 0 \hbox{ in } \mathbb R^+\times \mathbb R^3\\
    u(0,x) =u^0(x) \;\;\hbox{ in }\mathbb R^3,
  \end{cases}
$$
we refer the reader to \cite{HK}. Clearly, when $\alpha=0$ or
$\beta=1$ it is corresponds to the classical anisotropic
Navier-Stokes equation for more details the reader is referenced
to the book \cite{JP} . The second purpose of this paper is to
study the system $(NS_2)$ in the anisotropic Sobolev space
$H^{0,1}(\R^3)$. Before stating the second main result, we define
the space corresponding to the system:
$$\begin{array}{lcl}\mathcal G_{\beta}&=&\{f:\mathbb
R^+\times\mathbb R^3\rightarrow\mathbb R^3\,measurable;\\
&&(e^{\beta |f|^2}-1)|f|^2,e^{\beta
|f|^2}(\partial_{3}(|f|^2))^2,(e^{\beta |f|^2}-1)|\partial_{3}
f|^2\in L^1_{loc}(\R^+,L^1(\R^3))\}.\end{array}$$
\begin{theorem}\label{th2} Let $u^0\in H^{0,1}(\mathbb R^3)$ be a divergence free vector fields, then there is a unique global solution of $(NS_2)$:
$u\in L_{loc}^\infty(\R^+,H^{0,1}(\mathbb R^3)\cap
C(\R^+,L^2(\R^3))\cap L_{loc}^2(\R^+,H^{1,1}(\mathbb R^3))\cap
\mathcal G_{\beta}$. Moreover, for all $t\geq0$
\begin{equation}\label{eqth21}\|u(t)\|_{L^2}^2+2\int_0^t\|\nabla_h
u\|_{L^2}^2+2\alpha\int_0^t\|(e^{\beta |u|^2}-1)|u|^2\|_{L^1}\leq
\|u^0\|_{L^2}^2,\end{equation}
\begin{equation}\label{eqth22}\| \partial_{3} u(t)\|_{L^2}^2+\int_0^t\|\nabla_{h}\partial_{3} u\|_{L^2}^2dz+
\alpha\int_0^t\|(e^{\beta |u|^2}-1)|\partial_{3} u|^2\|_{L^1}~\quad\quad\quad\quad\quad\quad\quad\quad\quad~\end{equation}
$$~\quad\quad\quad\quad\quad\quad\quad\quad\quad~+\alpha \beta \int_0^t\|e^{\beta |u|^2}(\partial_{3}(|u|^2)^{2}\|_{L^1}\leq \|\partial_{3} u^0\|_{L^2}^2 e^{\frac{6t}{\alpha\beta^2}}.$$
\end{theorem}
\begin{rem}\label{rem2}\begin{enumerate}
\item Combining the above result $u\in C(\R^+,L^2(\R^3))$ and the fact $u\in L^\infty_{loc}(\mathbb R^+,H^{0,1}(\R^3))$ with the interpolation result, we get: For all $s<1$, we have $u\in C(\mathbb R^+,H^{0,s}(\R^3))$.
\item By inequalities (\ref{eqth21})-(\ref{eqth22}) and Proposition \ref{prop1}, we get
$$\lim_{t\rightarrow0}\|u(t)-u^0\|_{H^{0,1}}=0.$$
\item By using the same idea of Remark \ref{rem1}-(3), we get : $u\in C_r(\mathbb R^+,H^{0,1}(\R^3))$.
\item The continuity set of $u$ solution of $(NS_2)$ in $H^{0,1}(\R^3)$: Put the following subset of $\R^+$
$$A'=\{t\in\R^+;\;u ~~~\mbox{discontinous at} ~~~t~~~\mbox{in}\,H^{0,1}(\R^3)\}.$$
$A'$ is at most countable.
\end{enumerate}
\end{rem}

The remainder of our paper is organized as follows. In the second
section we give some notations, definitions and preliminary
results. Throughout Section 3, we will study the uniqueness and global existence of solution
of Cauchy problem $(NS_1)$ in $H^1(\R^3)$. In section 4, we will interest to anisotropic case :
We will also prove the uniqueness and global existence of solution in $H^{0,1}(\R^3)$.
\section{\bf Notations and preliminary results}
\subsection{Notations} In this section, we collect some notations and definitions that will be used later.\\
\begin{enumerate}
\item[$\bullet$] The Fourier transformation is normalized as
$$
\mathcal{F}(f)(\xi)=\widehat{f}(\xi)=\int_{\mathbb R^3}\exp(-ix.\xi)f(x)dx,\,\,\,\xi=(\xi_1,\xi_2,\xi_3)\in\mathbb R^3.
$$
\item[$\bullet$] The inverse Fourier formula is
$$
\mathcal{F}^{-1}(g)(x)=(2\pi)^{-3}\int_{\mathbb R^3}\exp(i\xi.x)g(\xi)d\xi,\,\,\,x=(x_1,x_2,x_3)\in\mathbb R^3.
$$
\item[$\bullet$] The convolution product of a suitable pair of function $f$ and $g$ on $\mathbb R^3$ is given by
$$
(f\ast g)(x):=\int_{\mathbb R^3}f(y)g(x-y)dy.
$$
\item[$\bullet$] If $f=(f_1,f_2,f_3)$ and $g=(g_1,g_2,g_3)$ are two vector fields, we set
$$
f\otimes g:=(g_1f,g_2f,g_3f),
$$
and
$$
{\rm div}\,(f\otimes g):=({\rm div}\,(g_1f),{\rm div}\,(g_2f),{\rm div}\,(g_3f)).
$$
Moreover, if $\rm{div}\,g=0$ we obtain
$$
{\rm div}\,(f\otimes g):=g_1\partial_1f+g_2\partial_2f+g_3\partial_3f:=g.\nabla f.
$$
\item[$\bullet$] Let $(B,||.||)$, be a Banach space, $1\leq p \leq\infty$ and  $T>0$. We define $L^p_T(B)$ the space of all
measurable functions $[0,t]\ni t\mapsto f(t) \in B$ such that $t\mapsto||f(t)||\in L^p([0,T])$.\\
\item[$\bullet$] The Sobolev space $H^s(\mathbb R^3)=\{f\in \mathcal S'(\mathbb R^3);\;(1+|\xi|^2)^{s/2}\widehat{f}\in L^2(\mathbb R^3)\}$.\\
\item[$\bullet$] The homogeneous Sobolev space $\dot H^s(\mathbb R^3)=\{f\in \mathcal S'(\mathbb R^3);\;\widehat{f}\in L^1_{loc}\;{\rm and}\;|\xi|^s\widehat{f}\in L^2(\mathbb R^3)\}$.\\
\item[$\bullet$] For $R>0$, the Friedrich operator $J_R$ is defined by
$$J_R(D)f=\mathcal F^{-1}({\bf 1}_{|\xi|<R}\widehat{f}).$$
\item[$\bullet$] The Leray projector $\mathbb P:(L^2(\R^3))^3\rightarrow (L^2(\R^3))^3$ is defined by
$$\mathcal F(\mathbb P f)=\widehat{f}(\xi)-(\widehat{f}(\xi).\frac{\xi}{|\xi|})\frac{\xi}{|\xi|}=M(\xi)\widehat{f}(\xi);\;M(\xi)=(\delta_{k,l}-\frac{\xi_k\xi_l}{|\xi|^2})_{1\leq k,l\leq 3}.$$
\item[$\bullet$] $L^2_\sigma(\R^3)=\{f\in (L^2(\R^3))^3;\;{\rm
div}\,f=0\}$. \item[$\bullet$] For $s_1,s_2\in\R$, the anisotropic
Sobolev spaces are defined by: $$H^{s_1,s_2}(\R^3)=\{f\in
S'(\R^3);\;(1+|\xi_h|^2)^{s_1/2}(1+\xi_3^2)^{s_2/2}\widehat{f}(\xi)\in
L^2(\R^3)\},$$
$$\dot H^{s_1,s_2}(\R^3)=\{f\in
S'(\R^3);\;|\xi_h|^{s_1}|\xi_3|^{s_2}\widehat{f}(\xi)\in
L^2(\R^3)\}.$$
\item[$\bullet$] Let $(B,||.||)$, be a Banach space
and  $I$ be nonempty interval. We define $C_{r}(I,B)$ the space of
all right continuous functions : $I\ni t\mapsto f(t) \in B$.
\item[$\bullet$] We often use the convex inequality: For $p,q\in(1,\infty)$ such that $\frac{1}{p}+\frac{1}{q}=1$, we have
$$ab\leq \frac{a^p}{p}+\frac{b^q}{q},\;\forall a,b\in\R^+.$$
\end{enumerate}
\subsection{Preliminaries}
In this section, we recall some classical results and we give new technical lemmas.
\begin{prop}(\cite{HBAF})\label{prop1} Let $H$ be a Hilbert space.
\begin{enumerate}
\item If $(x_n)$ is a bounded sequence of elements in $H$, then there is a subsequence $(x_{\varphi(n)})$ such that
$$(x_{\varphi(n)}|y)\rightarrow (x|y),\;\forall y\in H.$$
\item If $x\in H$ and $(x_n)$ is a bounded sequence of elements in $H$ such that
$$(x_n|y)\rightarrow (x|y),\;\forall y\in H.$$
Then $\|x\|\leq\displaystyle\liminf_{n\rightarrow\infty}\|x_n\|.$
\item If $x\in H$ and $(x_n)$ is a bounded sequence of elements in $H$ such that
$$\begin{array}{l}
(x_n|y)\rightarrow (x|y),\;\forall y\in H\\
\displaystyle\limsup_{n\rightarrow\infty}\|x_n\|\leq \|x\|,\end{array}$$
then $\displaystyle\lim_{n\rightarrow\infty}\|x_n-x\|=0.$
\end{enumerate}
\end{prop}
\begin{rem}\label{remprop1} Combining Proposition \ref{prop1}-(1) and (2), we get: If $(x_n)$ is a bounded sequence of elements in $H$ such that
$$(x_n|y)\rightarrow (x|y),\;\forall y\in H,$$
then
$$\lim_{n\rightarrow\infty}\|x_n-x\|=0\Longleftrightarrow \lim_{n\rightarrow\infty}\|x_n\|=\|x\|.$$
\end{rem}
\begin{lem}(\cite{JYC})\label{LP}
Let $s_1,\ s_2$ be two real numbers and $d\in\N$.
\begin{enumerate}
\item If $s_1<d/2$\; and\; $s_1+s_2>0$, there exists a constant  $C_1=C_1(d,s_1,s_2)$, such that: if $f,g\in \dot{H}^{s_1}(\mathbb{R}^d)\cap \dot{H}^{s_2}(\mathbb{R}^d)$, then $f.g \in \dot{H}^{s_1+s_2-\frac{d}{2}}(\mathbb{R}^d)$ and
$$\|fg\|_{\dot{H}^{s_1+s_2-\frac{d}{2}}}\leq C_1 (\|f\|_{\dot{H}^{s_1}}\|g\|_{\dot{H}^{s_2}}+\|f\|_{\dot{H}^{s_2}}\|g\|_{\dot{H}^{s_1}}).$$
\item If $s_1,s_2<d/2$\; and\; $s_1+s_2>0$ there exists a constant $C_2=C_2(d,s_1,s_2)$ such that: if $f \in \dot{H}^{s_1}(\mathbb{R}^d)$\; and\; $g\in\dot{H}^{s_2}(\mathbb{R}^d)$, then  $f.g \in \dot{H}^{s_1+s_2-1}(\mathbb{R}^d)$ and
$$\|fg\|_{\dot{H}^{s_1+s_2-\frac{d}{2}}}\leq C_2 \|f\|_{\dot{H}^{s_1}}\|g\|_{\dot{H}^{s_2}}.$$
\end{enumerate}
 \end{lem}
\begin{lem}\label{lem4}
    Let $\alpha>0$ and $d\in\N$. Then, for all $x,y\in\R^d$, we have
    $$  (|x|^{\alpha}x-|y|^{\alpha}y).(x-y)\geq C_{\alpha}(|x|^{\alpha}+|y|^{\alpha})|x-y|^{2}, $$
with $C_{\alpha}=\min(\frac{1}{18},\frac{1}{2^{\alpha+1}})>0$.
\end{lem}
{\bf Proof.} In all the proof we suppose that $|x|\geq |y|$. Particularly, we have
    $$|x|^{\alpha}=\frac{|x|^{\alpha}+|x|^{\alpha}}{2}\geq\frac{|y|^{\alpha}+|x|^{\alpha}}{2}.$$
{\bf First case:} we suppose that $(x,y)$ is {\bf related}. We treat two subcases:\\
$\bullet$ Suppose that $x.y\leq 0$, then $x.y=-|x|.|y|$ and $|x-y|=|x|+|y|\leq2|x|$. We have
        \begin{align*}
        (|x|^{\alpha}x-|y|^{\alpha}y).(x-y)&=|x|^{\alpha+2}+|y|^{\alpha+2}-(|x|^\alpha+|y|^\alpha)x.y\\
&=|x|^{\alpha+2}+|y|^{\alpha+2}+(|x|^\alpha+|y|^\alpha)|x|.|y|\\
&=(|x|^{\alpha+1}+|y|^{\alpha+1})(|x|+|y|)\\
&=(|x|^{\alpha+1}+|y|^{\alpha+1})|x-y|\\
&\geq |x|^{\alpha+1}|x-y|\\
&\geq |x|^{\alpha} |x||x-y|\\
&\geq\frac{|y|^{\alpha}+|x|^{\alpha}}{2}.\frac{|x-y|^{2}}{2}\\
&\geq\frac{1}{4}(|x|^{\alpha}+|y|^{\alpha}).|x-y|^{2}.
    \end{align*}
$\bullet$ Suppose that $xy>0$, then $x.y=|x|.|y|$ and $|x-y|=|x|-|y|$. We have
    \begin{align*}
    (|x|^{\alpha}x-|y|^{\alpha}y).(y-x)&=(|x|^{\alpha+1}-|y|^{\alpha+1}).|x-y|\\
&=(|x|^{\alpha+1}-|y|^{\alpha+1}).(|x|-|y|)\\
&=|x|^{\alpha+2}(1-(\frac{|y|}{|x|})^{\alpha+1})(1-\frac{|y|}{|x|}).
\end{align*} Let $\theta=\displaystyle\frac{|x|}{|y|}\in[0,1]$, then
\begin{align*}
    (|x|^{\alpha}x-|y|^{\alpha}y).(y-x)&\geq |x|^{\alpha+2}(1-\theta^{\alpha+1})(1-\theta)\\
&\geq |x|^{\alpha+2}(1-\theta)^2\\
&\geq |x|^{\alpha}(|x|-|x|\theta)^2\\
&\geq |x|^{\alpha}(|x|-|y|)^2\\
&\geq |x|^{\alpha}|x-y|^2\\
&\geq \frac{1}{2}(|x|^\alpha +|y|^{\alpha})|y-x|^2.
\end{align*}
{\bf Second case:} Suppose that $(x,y)$ are {\bf linearly independent} elements. There are two subcases:\\
$\bullet$ If $x.y\leq 0,$ then
    $$|x-y|^{2}=|x|^2+|y|^2-2x.y\Rightarrow|x-y|^{2}\leq 2(|x|^2+|y|^2)\leq 4|x|^2.$$
We have
    \begin{align*}
        (|x|^{\alpha}x-|y|^{\alpha}y).(x-y)&=|x|^{\alpha+2}+|y|^{\alpha+2}-|x|^{\alpha}x.y-|y|^{\alpha}x.y\\&\geq |x|^{\alpha+2}\\&\geq |x|^{\alpha} |x|^{2}\\&\geq\frac{|x|^{\alpha}+|x|^{\alpha}}{2}.\frac{|x-y|^{2}}{4}\\&\geq\frac{1}{8}(|x|^{\alpha}+|y|^{\alpha}).|x-y|^{2}.
    \end{align*}
$\bullet$ If $x.y>0,$ and  suppose that $\displaystyle|y|\leq \frac{|x|}{2}.$ We have
    $$|x-y|\leq |x|+|y|\leq \frac{3}{2}|x|\Rightarrow\frac{4}{9}|x-y|^{2}\leq|x|^2.$$
Put the following vectors:
$$v=\frac{x}{|x|},\;w=\frac{y}{|x|}.$$
Clearly, we get
$$|v|=1,\;|w|=\frac{|y|}{|x|}\in[0,\frac{1}{2}],\;v.w>0.$$
We have
    \begin{align*}
    (|x|^{\alpha}x-|y|^{\alpha}y).(x-y)&=|x|^{\alpha+2}(v-|w|^{\alpha}w).(v-w)\\
&=|x|^{\alpha+2}(1+|w|^{\alpha+2}-(1+|w|^\alpha)v.w)\\
&\geq|x|^{\alpha+2}(1+|w|^{\alpha+2}-(1+|w|)|v|.|w|)\\
&\geq|x|^{\alpha+2}(1+|w|^{\alpha+2}-\frac{3}{2}|v|.|w|)\\
&\geq|x|^{\alpha+2}(1+|w|^{\alpha+2}-\frac{3}{2}|w|)\\
&\geq|x|^{\alpha+2}(1+|w|^{\alpha+2}-\frac{3}{4})\\
&\geq|x|^{\alpha+2}(1-\frac{3}{4})\\
&\geq\frac{1}{4}|x|^{\alpha}|x|^2\\
&\geq\frac{1}{4}\frac{|x|^{\alpha}+|y|^\alpha}{2}\frac{4}{9}|x-y|^{2}\\
&\geq\frac{1}{18}(|x|^{\alpha}+|y|^\alpha)|x-y|^{2}.
 \end{align*}
$\bullet$ If $x.y>0,$ and  suppose that $\displaystyle\frac{|x|}{2}<|y|\leq|x|$.\\
Put the plan $P=Span\{x,y\}$ and
$\displaystyle B=(u_1=\frac{x}{|x|},u_2)$ a normalized basis of $P$.
We start by noting the following relations
    $$\begin{array}{l}
x=|x|u,\\
y=a_{1}u_1+a_{2}u_2,\\
a_{1}=x.y>0,\\
|y|=\sqrt{a_{1}^{2}+a_{2}^{2}}.
\end{array}$$
Then, if we put $\displaystyle z=\frac{y}{|x|}=b_{1}u_1+b_{2}u_2$,
we obtain
$$\begin{array}{l}
\displaystyle b_{1}=\frac{a_1}{|x|}>0,\\
\displaystyle\frac{1}{2}\leq |z|=\sqrt{b_{1}^{2}+b_{2}^{2}}\leq1.
\end{array}$$
Then
    \begin{align*}
         (|x|^{\alpha}x-|y|^{\alpha}y).(x-y)
&=|x|^{\alpha+2}[(u-|z|^{\alpha}z).(u-z)]\\&=|x|^{\alpha+2}[(u-|z|^{\alpha}(b_{1}u+b_{2}v)).(u-(b_{1}u+b_{2}v))]\\
&=|x|^{\alpha+2}[((1-|z|^{\alpha}b_{1})u-|z|^{\alpha}b_{2}v).((1-b_{1})u-b_{2}v)]\\
&=|x|^{\alpha+2}[(1-|z|^{\alpha}b_{1})(1-b_{1})+|z|^{\alpha}b_{2}^{2}]\\
&\geq |x|^{\alpha+2}[(1-b_{1})^2+(\frac{1}{2})^{\alpha}b_{2}^{2}],\;(b_{1}>0\,\mbox{and}\,\frac{1}{2}\leq |z|\leq 1)\\
&\geq (\frac{1}{2})^{\alpha}|x|^{\alpha+2}[(1-b_{1})^2+b_{2}^{2}]\\
&\geq(\frac{1}{2})^{\alpha}|x|^{\alpha+2}|u-z|^2\\
&\geq(\frac{1}{2})^{\alpha}|x|^{\alpha}\Big||x|u-|x|z\Big|^2\\
&\geq(\frac{1}{2})^{\alpha}|x|^{\alpha}|x-y|^2\\
&\geq \frac{1}{2}(\frac{1}{2})^{\alpha}(|x|^{\alpha}+|y|^{\alpha})|x-y|^2\\
         &\geq (\frac{1}{2})^{\alpha+1}(|x|^{\alpha}+|y|^{\alpha})|x-y|^2.
    \end{align*}
So, the real $C_{\alpha}=\min(\frac{1}{18},b_{\alpha}
=\frac{1}{2^{\alpha+1}})$ answers the question.
 \begin{rem}\label{rem11} It's easy to see that $b_{\alpha}$ is decreasing to $0$. By a
 straightforward computation we get, for $\alpha\geq 3$,  $\displaystyle b_{\alpha}<\frac{1}{18}$.
If $k\in\N$, then $C_{2k}=\left\{\begin{array}{l}
\displaystyle \frac{1}{18},\;if\;k=1\\
\displaystyle b_{2k}=\frac{1}{2^{2k+1}},\;if \;k\geq2.
\end{array}\right.$
\end{rem}
Then
$$\forall k\in\N;\;C_{2k}\geq \frac{8}{18}b_{2k}=\frac{2}{9}\frac{1}{2^{2k}}=\frac{2}{9}\frac{1}{4^{k}}.$$
Combining Lemma \ref{lem4} and Remark \ref{rem1}, we get the following result.
\begin{lem}\label{lem44}
If $\beta>0$, then, for all $x,y\in\R^d$, we have
$$\Big((e^{\beta|x|^2}-1)x-(e^{\beta|y|^2}-1)y\Big).(x-y)\geq \frac{2}{9}\Big[(e^{\frac{\beta}{4}|x|^2}-1)+(e^{\frac{\beta}{4}|y|^2}-1)\Big]|x-y|^2.$$
\end{lem}
In the following we give another version of Gronwall's lemma that we often use:
\begin{lem}\label{lem45}
Let $A,T>0$ and $f,g,h:[0,T]\rightarrow\R^+$ three continuous functions such that
\begin{equation}\label{LG2}\forall t\in[0,T];\;f(t)+\int_0^tg(z)dz\leq A+\int_0^th(z)f(z)dz.\end{equation}
Then $$\forall t\in[0,T];\;f(t)+\int_0^tg(z)dz\leq A\exp(\int_0^th(z)dz).$$
\end{lem}
{\bf Proof.} By Gronwall Lemma, we get
$$\forall t\in[0,T];\;f(t)\leq A\exp(\int_0^th(z)dz).$$
Put this inequality in $(\ref{LG2})$ we obtain
$$\begin{array}{lcl}
\displaystyle f(t)+\int_0^tg(z)dz&\leq&\displaystyle  A+\int_0^th(z)A\exp(\int_0^zh(r)dr)dz\\
&\leq&\displaystyle  A+A\int_0^th(z)\exp(\int_0^zh(r)dr)dz\\
&\leq&\displaystyle  A+A\int_0^t\Big(\exp(\int_0^zh(r)dr)\Big)'dz\\
&\leq&\displaystyle  A+A\Big(\exp(\int_0^th(r)dr)-1\Big)\\
&\leq&\displaystyle  A\exp(\int_0^th(r)dr),
\end{array}$$
which ends the proof.\\

To prove the right continuity of strong solutions of $(NS_1)$ and $(NS_2)$, I need the following classical result:
\begin{lem}\label{lem5}
Let $f:I\rightarrow\R$ be a monotonic function on an interval $I$.
Then there is $A\subset\R$ at most countable family such that $f$
is continuous on $I\setminus A$.
\end{lem}
{\bf Proof.} Suppose that $f$ is increasing (if $f$ is decreasing, we can consider
$g=-f$). Then, for $t\in int(I)$, we have $$f\;{\rm
is\,discontinuous\,at}\,t\Longleftrightarrow
\underset{t^-}{\lim}f<\underset{t^+}{\lim}f.$$ Let $A=\{t\in\R;\;
f~~\mbox{discontious at}~~t\}$ and $a\in A\cap int(I)$, then we
have $\underset{a^-}{\lim}f<\underset{a^+}{\lim}f$. So
$(\underset{a^-}{\lim}f,\underset{a^+}{\lim}f)\cap
\mathbb{Q}\neq\varnothing$, and we can choose
$r_{a}\in(\underset{a^-}{\lim}f,\underset{a^+}{\lim}f)\cap
\mathbb{Q}$. Then, the following function
    $$\left.
    \begin{array}{cc}
        \varphi :& A \rightarrow \mathbb{Q}\\
        & a \longmapsto r_{a}
    \end{array}\right. $$
    is well defined. For  $a,b\in A$ such that $a<b$, we have
$$r_{a}<\lim_{a^+}f\leq\lim_{b^-}f<r_b$$ so $\varphi(a)<\varphi(b)$
which implies that $\varphi$ is injective function. Therefore $A$
is at most countable family.
\section{\bf Proof of Theorem \ref{th1}}
\subsection{A priori estimates }
We start by taking the $L^2$ scalar product of the first equation of $(NS_1)$ with $u$, we get
\begin{equation}\label{ineqab1}
\|u(t)\|_{L^2}^2+2\int_0^t\|\nabla u\|_{L^2}^2+2\alpha\int_0^t\|(e^{\beta |u|^2}-1)|u|^2\|_{L^1}\leq \|u^0\|_{L^2}^2.
\end{equation}
Also, taking the $\dot{H}^{1}$ scalar product of the first equation of $(NS_1)$ with $u$, we obtain
$$\langle \partial_{t} \nabla u,\nabla u\rangle_{L^{2}}-  \langle \Delta \nabla u,\nabla u\rangle_{L^{2}}+\alpha \langle \nabla ((e^{\beta|u|^{2}}-1)u),\nabla u\rangle_{L^{2}}\leq |\langle u\nabla u,\Delta u\rangle_{L^{2}}|.$$
By using the following identities,
$$\begin{array}{ccl}
\alpha \langle \nabla ((e^{\beta|u|^{2}}-1)u),\nabla u\rangle_{L^{2}}&=&\displaystyle\frac{\alpha\beta}{2}\|e^{\beta|u|^{2}}|\nabla(|u|^{2})|^{2}\|_{L^{1}}
+\alpha\|(e^{\beta|u|^{2}}-1)|\nabla u|^{2}\|_{L^{1}}\\
\displaystyle|\langle u\nabla u,\Delta u\rangle_{L^{2}}|&\leq&\displaystyle\int_{\R^3}|u|.|\nabla u|.|\Delta u|\\
&\leq&\displaystyle\frac{1}{2}\int_{\R^3}|u|^2|\nabla u|^2+\frac{1}{2}\int_{\R^3}|\Delta u|^2,
\end{array}$$
we obtain
    $$\frac{d}{dt}\|\nabla u\|^{2}_{L^{2}}+\|\Delta u\|^{2}_{L^{2}}+\alpha \beta\|e^{\beta|u|^{2}}|\nabla(|u|^{2})|^{2}\|_{L^{1}}
+2\alpha \|(e^{\beta|u|^{2}}-1)|\nabla u|^{2}\|_{L^{1}}\leq\int_{\R^{3}}|u|^{2}|\nabla u|^{2}.
$$
Moreover, the elementary inequalities
$$\begin{array}{ccl}
    \alpha (e^{\beta|u|^{2}}-1)&\geq&\displaystyle\alpha\frac{(\beta|u|^{2})^2}{2!}\\
    &\geq&\displaystyle\frac{\alpha\beta^2}{2}|u|^{4}\\
    ~\quad\quad\quad\quad~|u|^{2}&=&\displaystyle(\frac{\sqrt{\alpha}\beta|u|^{2}}{2\sqrt{2}}).(\frac{\sqrt{2}}{\sqrt{\alpha}\beta})\\
    &\leq&\displaystyle\frac{\alpha\beta^{2}}{16}|u|^{4}+\frac{1}{\alpha\beta^{2}}
    \end{array}$$
imply
$$\frac{d}{dt}\|\nabla u\|^{2}_{L^{2}}+\|\Delta  u\|^{2}_{L^{2}}+\alpha \beta\|e^{\beta|u|^{2}}|\nabla(|u|^{2})|^{2}\|_{L^{1}}+\alpha \|(e^{\beta|u|^{2}}-1)|\nabla u|^{2}\|_{L^{1}}\leq \frac{1}{\alpha \beta^{2}}\|\nabla u\|^{2}_{L^{2}}.$$
Integrate on $[0,t]$, we get
\begin{align*}
    \|\nabla u\|^{2}_{L^{2}}+\int_{0}^{t}\|\Delta u\|^{2}_{L^{2}}+\alpha \beta\int_{0}^{t}\|e^{\beta|u|^{2}}|\nabla(|u|^{2})|^{2}\|_{L^{1}}\\+\int_{0}^{t}\alpha \|(e^{\beta|u|^{2}}-1)|\nabla u|^{2}\|_{L^{1}}&\leq \|\nabla u^{0}\|^{2}_{L^{2}}+\frac{1}{\alpha \beta^{2}}\int_{0}^{t}\|\nabla u\|^{2}_{L^{2}}.
    \end{align*}
This inequality implies two results, the first by applying Lemma \ref{lem45}:
\begin{equation}\label{ineqab2}
\|\nabla u\|^{2}_{L^{2}}+\int_{0}^{t}\|\Delta u\|^{2}_{L^{2}}+\alpha \beta\int_{0}^{t}\|e^{\beta|u|^{2}}|\nabla(|u|^{2})|^{2}\|_{L^{1}}~\quad\quad\quad\quad\quad\quad\quad\quad\quad\quad\quad\quad~\end{equation}$$~\quad\quad\quad\quad\quad\quad\quad\quad\quad\quad\quad\quad~+\int_{0}^{t}\alpha \|(e^{\beta|u|^{2}}-1)|\nabla u|^{2}\|_{L^{1}}\leq \|\nabla u^{0}\|^{2}_{L^{2}}e^{\frac{t}{\alpha\beta^{2}}},
$$
and the second by using inequality (\ref{ineqab1}):
\begin{equation}\label{ineqab3}\|\nabla u\|^{2}_{L^{2}}+\int_{0}^{t}\|\Delta u\|^{2}_{L^{2}}+\alpha \beta\int_{0}^{t}\|e^{\beta|u|^{2}}|\nabla(|u|^{2})|^{2}\|_{L^{1}}~\quad\quad\quad\quad\quad\quad\quad\quad\quad\quad\quad\quad~\end{equation}
$$~\quad\quad\quad\quad\quad\quad\quad\quad\quad\quad\quad\quad~+\alpha\int_{0}^{t} \|(e^{\beta|u|^{2}}-1)|\nabla u|^{2}\|_{L^{1}}\leq \|\nabla u^{0}\|^{2}_{L^{2}}+\frac{\|u^{0}\|^{2}_{L^{2}}}{\alpha \beta^{2}}.
$$
Absolutely, these bounds come from the approximate solutions via the Friedrich's regularization procedure. Hence it remains to pass to the limit in the sequence of solutions of approximate schema. The passage to the limit follows using classical argument by combining Ascoli's theorem and the Cantor diagonal process (see \cite{JB}). Moreover, inequalities (\ref{eqth1})-(\ref{eqth2})-(\ref{eqth3}) are given by (\ref{ineqab1})-(\ref{ineqab2})-(\ref{ineqab3}). Finally, $u$ is in the space $C(\R^+,H^{-2}(\R^3))$  as an interpolation  between $C(\R^+,H^{-4}(\R^3))$ and $L^\infty(\R^+,H^{1}(\R^3))$.
\subsection{Uniqueness}: Let $u,v$ two solutions of $(NS_1)$ and put $w=u-v$. We make the difference of two following equations
$$\begin{array}{lcl}\partial_t u
	-\Delta u+ u.\nabla u  +\alpha (e^{\beta |u|^2}-1)u&=&\;\;-\nabla p\\
	\partial_t v
	-\Delta v+ v.\nabla v +\alpha (e^{\beta |v|^2}-1)v&=&\;\;-\nabla \tilde p,
\end{array}
$$
we get
$$\partial_t w
-\Delta w+ w.\nabla u +v.\nabla w +\alpha (e^{\beta |u|^2}-1)u-\alpha(e^{\beta |v|^2}-1)v=-\nabla (p-\tilde p).$$
Taking the $L^2$ scalar product with $w$, we obtain :
$$\langle \partial_t w,w\rangle_{L^2}
-\langle \Delta w,w\rangle_{L^{2}} +\langle \alpha (e^{\beta |u|^2}-1)u-\alpha(e^{\beta |v|^2}-1)v,w\rangle_{L^2}\leq |\langle w.\nabla u,w\rangle_{L^2}|.$$
Using Lemma \ref{lem44}, we get
$$\begin{array}{lcl}
\displaystyle\frac{1}{2} \frac{d}{dt} \|w\|^{2}_{L^{2}}+\|\nabla w\|^{2}_{L^{2}}+\frac{2\alpha}{9}\int_{\R^{3}}(e^{\frac{\beta}{4}|u|^2}-1)|w|^{2}&\leq&\|w u\|_{L^{2}}\|\nabla w\|_{L^{2}}\\
&\leq&\displaystyle\frac{1}{2} \|w u\|^2_{L^{2}}+\frac{1}{2}\|\nabla w\|^2_{L^{2}}
\end{array}$$
and
$$\frac{1}{2} \frac{d}{dt} \|w\|^{2}_{L^{2}}+\frac{1}{2}\|\nabla w\|^{2}_{L^{2}}+\frac{2\alpha}{9}\int_{\R^{3}}(e^{\frac{\beta}{4}|u|^2}-1)|w|^{2}\leq \frac{1}{2} \int_{\R^{3}} |w|^2 |u|^2dz.$$
By the elementary inequalities
$$\begin{array}{ccl}
\displaystyle\frac{2\alpha}{9}(e^{\frac{\beta}{4}|u|^2}-1)&\geq&\displaystyle\frac{2\alpha}{9}.\frac{1}{2!}(\frac{\beta}{4}|u|^2)^2=\frac{\alpha\beta^2}{36}|u|^4\\
~\;\;\;\;\;\;\;\;\;\;\;\;\;\;\;~|u|^2&=&\displaystyle\Big(\frac{\sqrt{\alpha}\beta}{\sqrt{18}}|u|^2\Big).\Big(\frac{\sqrt{18}}{\sqrt{\alpha}\beta}\Big)\leq \frac{\alpha\beta^2}{36}|u|^4+\frac{9}{\alpha\beta^2},
\end{array}$$
we get
\begin{align*}
	 \frac{1}{2} \frac{d}{dt} \|w\|^{2}_{L^{2}}+\frac{1}{2}\|\nabla w\|^{2}_{L^{2}}\leq \frac{18}{\alpha\beta^2} \int_{\R^{3}} |w|^2.
\end{align*}	
According to Gronwall Lemma, we obtain
$$\|w(t)\|^{2}_{L^{2}}\leq \|w(0)\|^2_{L^{2}}e^{\frac{18}{\alpha\beta^2}t}.$$	
But $w(0)=0$, so $u=v$, which ends the proof.
\section{\bf Proof of Theorem \ref{th2}.}
In this section, we do the same procedure as the previous proof.
\subsection{A priori estimates} We start by taking the scalar product in $L^2(\R^3)$
\begin{equation}\label{ineqac1}
\|u(t)\|_{L^2}^2+2\int_0^t\|\nabla_h u\|_{L^2}^2+2\alpha\int_0^t\|(e^{\beta |u|^2}-1)|u|^2\|_{L^1}\leq \|u^0\|_{L^2}^2.
\end{equation}
At present, taking $\dot H^{0,1}$ scalar product of the system with the solution $u$, we get:
$$\frac{1}{2}\frac{d}{dt}\|\partial_3 u\|^{2}_{L^{2}}+\|\nabla_{h}\partial_3 u\|^{2}_{L^{2}}+\frac{\alpha\beta}{2}\|e^{\beta |u|^{2}}|\partial_3 (|u|^{2})|^2\|_{L^{1}}~\quad\quad\quad\quad\quad\quad\quad\quad\quad\quad~$$$$~\quad\quad\quad\quad\quad\quad\quad\quad\quad\quad~+\alpha \|(e^{\beta|u|^{2}}-1)|\partial_3 u|^{2}\|_{L^{1}}=\langle \partial_3 u\nabla u,\partial_3 u\rangle_{L^{2}}.
$$
By following the process of \cite{HB}(see pages 1819-1820), we get two universal constants $C_0,C_1$ and
$$|\langle \partial_3 u\nabla u,\partial_3 u\rangle_{L^{2}}|\leq \frac{1}{2}\|\nabla_{h}\partial_{3}u\|_{L^{2}}^2+C_0\|\partial_{3}u\|_{L^{2}}^2 + C_1\||u|^2\partial_{3}u\|^{2}_{L^{2}}.
$$
Then, we obtain
$$
    \frac{1}{2}\frac{d}{dt}\|\partial_3 u\|^{2}_{L^{2}}+\|\nabla_{h}\partial_3 u\|^{2}_{L^{2}}+\frac{\alpha\beta}{2}\|e^{\beta |u|^{2}}|\partial_3 (|u|^{2})|\|_{L^{1}}~\quad\quad\quad\quad\quad\quad\quad\quad\quad\quad~$$
    $$~\quad\quad\quad\quad\quad\quad\quad\quad\quad\quad~+\alpha \|(e^{\beta|u|^{2}}-1)|\partial_3 u|^{2}\|_{L^{1}}\leq C_0\|\partial_{3}u\|_{L^{2}}^2 + C_1\||u|^2\partial_{3}u\|^{2}_{L^{2}}.
$$
Moreover, the elementary inequalities
$$\begin{array}{ccl}
    \displaystyle\frac{\alpha}{2} (e^{\beta|u|^{2}}-1)&\geq&\displaystyle\frac{\alpha}{2}\frac{(\beta|u|^{2})^3}{3!}\\
    &\geq&\displaystyle\frac{\alpha\beta^3}{12}|u|^{6}\\
    ~\quad\quad\quad\quad~C_1|u|^{2}&=&\displaystyle(\frac{\alpha^{1/3}\beta|u|^{2}}{4^{1/3}}).(\frac{C_14^{1/3}}{\alpha^{1/3}\beta})\\
    &\leq&\displaystyle\frac{\alpha\beta^{3}}{12}|u|^{6}+C(\alpha,\beta),\;\;C(\alpha,\beta)=\frac{256C_1^{4/3}}{\alpha^4\beta^{4/3}},
    \end{array}$$
imply
$$
\frac{1}{2}\frac{d}{dt}\|\partial_3 u\|^{2}_{L^{2}}+\|\nabla_{h}\partial_3 u\|^{2}_{L^{2}}+\frac{\alpha\beta}{2}\|e^{\beta |u|^{2}}|\partial_3 (|u|^{2})|\|_{L^{1}}~\quad\quad\quad\quad\quad\quad\quad\quad\quad\quad~$$$$~\quad\quad\quad\quad\quad\quad\quad\quad\quad\quad~+\frac{\alpha}{2} \|(e^{\beta|u|^{2}}-1)|\partial_3 u|^{2}\|_{L^{1}}\leq a_{\alpha,\beta}\|\partial_3 u\|_{L^2}^{2},
$$
where $a_{\alpha,\beta}=2(C_0+C(\alpha,\beta))$.
Integrate on $[0,t]$ and using Lemma \ref{lem45}, we get
$$
    \|\partial_3 u\|^{2}_{L^{2}}+\int_{0}^{t}\|\nabla_{h}\partial_{3}u\|^{2}_{L^{2}}+\alpha\beta\int_{0}^{t}\|e^{\beta |u|^{2}}|\partial_3 (|u|^{2})|\|_{L^{1}} ~\quad\quad\quad\quad\quad\quad\quad\quad\quad\quad\quad\quad~$$$$~\quad\quad\quad\quad\quad\quad\quad\quad\quad\quad\quad~+\alpha\int_{0}^{t} \|(e^{\beta|u|^{2}}-1)|\partial_3 u|^{2}\|_{L^{1}}\leq \|\partial_3 u^{0}\|^{2}_{L^{2}}\exp(a_{\alpha,\beta}t).
$$
\subsection{Approximate system} In this step we construct a global solution of $(NS_2)$, where we use a method inspired by \cite{JB}-\cite{JYC}. For this, consider the approximate system with the parameter $n\in\N$:
$$(S_n)
  \begin{cases}
     \partial_t u
 -\Delta_h J_nu+ J_n(J_nu.\nabla J_nu  +\alpha J_n[(e^{\beta |J_nu|^2}-1)J_nu] =\;\;-\nabla p_n\hbox{ in } \mathbb R^+\times \mathbb R^3\\
 p_n=(-\Delta)^{-1}\Big({\rm div}\,J_n(J_nu.\nabla J_nu  +\alpha {\rm div}\,J_n[(e^{\beta |J_nu|^2}-1)J_nu]\Big)\\
     {\rm div}\, u = 0 \hbox{ in } \mathbb R^+\times \mathbb R^3\\
    u(0,x) =J_nu^0(x) \;\;\hbox{ in }\mathbb R^3.
  \end{cases}
$$
By Cauchy-Lipschitz Theorem, we obtain a unique solution $u_n\in C^1(\R^+,L^2_\sigma(\R^3))$ of $(S_n)$. Moreover, $J_nu_n=u_n$. By the last section, we get
\begin{equation}\label{thapp1}\|u_n(t)\|_{L^2}^2+2\int_0^t\|\nabla_h u_n\|_{L^2}^2+2\alpha\int_0^t\|(e^{\beta |u_n|^2}-1)|u_n|^2\|_{L^1}\leq \|u^0\|_{L^2}^2.\end{equation}
\begin{equation}\label{thapp2}
 \|\partial_3 u_n(t)\|^{2}_{L^{2}}+\int_{0}^{t}\|\nabla_{h}\partial_{3}u_n\|^{2}_{L^{2}}+\alpha\beta\int_{0}^{t}\|e^{\beta |u_n|^{2}}|\partial_3 (|u_n|^{2})|\|_{L^{1}} ~\quad\quad\quad\quad\quad\quad\quad~\end{equation}$$~\quad\quad\quad\quad\quad\quad\quad\quad\quad\quad\quad~+\alpha\int_{0}^{t} \|(e^{\beta|u_n|^{2}}-1)|\partial_3 u_n|^{2}\|_{L^{1}}\leq \|\partial_3 u^{0}\|^{2}_{L^{2}}\exp(a_{\alpha,\beta}t).$$
Let $(T_q)_q\in(0,\infty)^\N$ such that $T_q<T_{q+1}$ and $T_q\rightarrow\infty$ as $q\rightarrow\infty$. Let $(\theta_q)_{q\in\N}$ be a sequence in $C_0^\infty(\R^3)$ such that: for all $q\in\N$
$$\left\{\begin{array}{l}
\theta_q(x)=1,\;\forall x\in B(0,q+1+\frac{1}{4})\\
\theta_q(x)=0,\;\forall x\in B(0,q+2)^c\\
0\leq \theta_q\leq 1.
\end{array}\right.$$
Using (\ref{thapp1})-(\ref{thapp2}) and classical argument by combining Ascoli's theorem and the Cantor diagonal process, we get a nondecreasing $\varphi:\N\rightarrow\N$ and
$u\in L^\infty(\R^+,H^{0,1}(\R^3))\cap C(\R^+,H^{-3}(\R^3))$ such that: for all $q\in\N$, we have
\begin{equation}\label{theq6}\lim_{n\rightarrow\infty}\|\theta_q(u_{\varphi(n)}-u)\|_{L^\infty([0,T_q],H^{-4})}=0,\end{equation}
and
\begin{equation}\label{thapp11}\|u(t)\|_{L^2}^2+2\int_0^t\|\nabla_h u\|_{L^2}^2+2\alpha\int_0^t\|(e^{\beta |u|^2}-1)|u|^2\|_{L^1}\leq \|u^0\|_{L^2}^2.\end{equation}

\begin{equation}\label{thapp12}
 \|\partial_3 u(t)\|^{2}_{L^{2}}+\int_{0}^{t}\|\nabla_{h}\partial_{3}u\|^{2}_{L^{2}}+\alpha\beta\int_{0}^{t}\|e^{\beta |u|^{2}}|\partial_3 (|u|^{2})|\|_{L^{1}} ~\quad\quad\quad\quad\quad\quad\quad~\end{equation}$$~\quad\quad\quad\quad\quad\quad\quad\quad\quad\quad\quad~+\alpha\int_{0}^{t} \|(e^{\beta|u|^{2}}-1)|\partial_3 u|^{2}\|_{L^{1}}\leq \|\partial_3 u^{0}\|^{2}_{L^{2}}\exp(a_{\alpha,\beta}t).$$
\noindent$\bullet$ We have $u\in C(\R^+,L^2(\R^3))$ (for the proof see Appendix). Then, by using interpolation inequality, we get $u\in C(\R^+,H^{0,s}(\R^3))$, for all $s<1$.
\subsection{Uniqueness} This proof is inspired by \cite{HB}. Let $u,\,v$ two solutions of $(NS_2)$ and denote by $w=u-v$. We make the difference of the following equations
$$
\begin{array}{lcl}
    \partial_t u
    -\Delta_{h} u+ u.\nabla u  +\alpha (e^{\beta |u|^2}-1)u &=&\;\;-\nabla p\\
    \partial_t v
    -\Delta_{h} v+ v.\nabla v +\alpha (e^{\beta |v|^2}-1)v &=&\;\;-\nabla \tilde p
\end{array}
$$
we get
$$\partial_t w
-\Delta_{h} w+ w.\nabla u +v.\nabla w +\alpha (e^{\beta |u|^2}-1)u-\alpha(e^{\beta |v|^2}-1)v=-\nabla (p-\tilde p).$$
Taking the $L^2$ scalar product with $w$, we obtain
$$\langle \partial_t w,w\rangle_{L^2}
-\langle \Delta_{h} w,w\rangle_{L^{2}} +\alpha\langle  (e^{\beta |u|^2}-1)u-(e^{\beta |v|^2}-1)v,w\rangle_{L^2}\leq |\langle w.\nabla u,w\rangle_{L^2}|$$
Using Lemma \ref{lem44}, we get:
\begin{align*}
    \frac{1}{2} \frac{d}{dt} \|w\|^{2}_{L^{2}}+\|\nabla_{h} w\|^{2}_{L^{2}}+\frac{2\alpha}{9}\int_{\R^{3}}\Big((e^{\frac{\beta}{4}|u|^2}-1)+(e^{\frac{\beta}{4}|v|^2}-1)\Big)|w|^{2}&\leq |\int_{\R^3} (w.\nabla u).w dx|
\end{align*}
and
\begin{align*}
    \frac{1}{2} \frac{d}{dt} \|w\|^{2}_{L^{2}}+\|\nabla_{h} w\|^{2}_{L^{2}}&\leq |\int_{\R^3} (w.\nabla u).w dx|.
\end{align*}
The second member of this inequality can be written as follows $$\int_{\R^3} (w.\nabla u).w dx=F_{1}+F_{2}$$
with $$\begin{array}{lcl}
F_1&=&\displaystyle\sum_{i=1}^{2}\sum_{j=1}^{3}\int_{\R^3}w_{i}\partial_{i}u_{j}w_{j}\\
F_2&=&\displaystyle\sum_{j=1}^{3}\int_{\R^3}w_{3}\partial_{3}u_{j}w_{j}.
\end{array}$$
By H\"older inequality, we have
$$F_{1}\leq \sum_{i=1}^{2}\sum_{j=1}^{3} \|\partial_{i}u_{j}\|_{L_{v}^\infty(L_{h}^2)}\|w_{i}\|_{L_{v}^2(L_{h}^4)}\|w_{j}\|_{L_{v}^2(L_{h}^4)}$$
Since $\dot{H}^\frac{1}{2}(\R^2)\hookrightarrow L^4(\R^2)$ and by interpolation we have :
$$\|w_{j}(.,x_3)\|_{L_{h}^4(\R^2)}\leq\|w_{j}(.,x_3)\|_{L^2(\R^2)}^{\frac{1}{2}}\|\nabla_h w_{j}(.,x_3)\|_{L^2(\R^2)}^{\frac{1}{2}}$$ 
so
$$\|w_{j}\|_{L^2(L_{h}^4)}\leq\|w_{j}\|_{L^2}^{\frac{1}{2}}\|\nabla_h w_{j}\|_{L^2}^{\frac{1}{2}}.$$
To estimate $\|\partial_{i}u_{j}\|_{L_{v}^\infty(L_{h}^2)}$, we write
$$\begin{array}{lcl}\|\partial_{i}u_{j}(.,x_3)\|^2_{L^2_h(\R^2)}&=&\displaystyle\int_{-\infty}^{x_3}\frac{d}{dz}\|\partial_{i}u_{j}(.,z)\|^2_{L^2_h(\R^2)}dz\\
&=&\displaystyle2\int_{-\infty}^{x_3}\langle\partial_{z} \partial_{i}u_{j}(.,z),\partial_{i}u_{j}(.,z)\rangle_{L^2_h(\R^2)} dz\\
&\leq&\displaystyle\|\partial_{3} \partial_{i}u_{j}\|_{L^2}\|\partial_{i}u_{j}\|_{L^2}.\end{array}$$
Then $$F_{1}\leq c\|\partial_{3}
\nabla_{h}u\|^{\frac{1}{2}}_{L^2}\|\nabla_{h}u\|^{\frac{1}{2}}_{L^2}\|w\|_{L^2}|\nabla_{h}
w\|_{L^2}.$$ By Young inequality, we obtain
\begin{align}\label{f1}F_{1}\leq \frac{1}{4}\|\nabla_{h} w\|^2_{L^2}+\frac{c}{4}(\|\partial_{3} \nabla_{h}u\|^{2}_{L^2}+\|\nabla_{h}u\|^{2}_{L^2})\|w\|^2_{L^2} \end{align}
The same procedure for $F_{2}$ we have :
\begin{align*}
    F_{2}\leq \|w_3\|_{L_{v}^\infty(L_h^2)}\|\nabla_{h} w\|^{\frac{1}{2}}_{L^2}\|\partial_{3} \nabla_{h}u\|^{\frac{1}{2}}_{L^2}\|\nabla_{h}u\|_{L^2}\|w\|_{L^2}.
    \end{align*}
Since
$$ \|w_3(.,x_3)\|_{L_h^2(\R^2)}^2=2\int_{-\infty}^{x_3}\int_{\R^2}w_3(x_h,z)\partial_{3}w_3(x_h,z)dx_hdz.$$ Using the fact that $\nabla.w=0$ so ${\rm div}_h\, w_h=-\partial_{3}w_3$ and 
$$ \begin{array}{lcl}
\|w_3(.,x_3)\|_{L_h^2(\R^2)}^2&=&\displaystyle-2\int_{-\infty}^{x_3}\int_{\R^2} w_3(x_h,z){\rm div_h} w_h (x_h,z)dx_h dz\\
&\leq&2\|{\rm div_h} w_h\|_{L^2(\R^3)}\|w_3\|_{L^2(\R^3)}.\end{array}$$
By Young inequality:
\begin{align}\label{f2}
    F_{2}\leq \frac{1}{4}\|\nabla_{h} w\|^2_{L^2}+\frac{c}{4}(\|\partial_{3} \nabla_{h}u\|^{2}_{L^2}+\|\partial_{3}u\|^{2}_{L^2})\|w\|^2_{L^2}
\end{align}
Hence, according to ($\ref{f1}$) and (\ref{f2}):
\begin{align*}
    \frac{1}{2} \frac{d}{dt} \|w\|^{2}_{L^{2}}+\frac{1}{2}\|\nabla_{h} w\|^{2}_{L^{2}}&\leq c(\|\partial_{3} \nabla_{h}u\|^{2}_{L^2}+\|\partial_{3}u\|^{2}_{L^2}+\|\nabla_{h}u\|^{2}_{L^2})\|w\|^2_{L^2}
\end{align*}
Integrate on $[0,t]$, we have :
\begin{align*}
\|w(t)\|^{2}_{L^{2}}+\int_{0}^{t}\|\nabla_{h} w\|^{2}_{L^{2}}&\leq\|w(0)\|^{2}_{L^{2}} +c\int_{0}^{t}(\|\partial_{3} \nabla_{h}u\|^{2}_{L^2}+\|\partial_{3}u\|^{2}_{L^2}+\|\nabla_{h}u\|^{2}_{L^2})\|w\|^2_{L^2}
\end{align*}
Now, by Gronwall Lemma :
$$\|w(t)\|^{2}_{L^{2}}\leq\|w(0)\|^{2}_{L^{2}} \exp\Big(c\int_{0}^{t}(\|\partial_{3} \nabla_{h}u\|^{2}_{L^2}+\|\partial_{3}u\|^{2}_{L^2}+\|\nabla_{h}u\|^{2}_{L^2})\Big).$$
Using inequalities (\ref{eqth21}) and (\ref{eqth22}), we get
\begin{equation}\label{eqf1}\|w(t)\|^{2}_{L^{2}}\leq\|w(0)\|^{2}_{L^{2}} \exp\Big(2c\|\partial_{3}u^0\|^{2}_{L^2}e^{\frac{6t}{\alpha\beta^2}}+c\|u^0\|^{2}_{L^2}\Big).\end{equation}
But $w(0)=u(0)-v(0)$, then $u=v$, which ends the proof of Theorem \ref{th2}.
\section{Appendix} In this section, we give a simple proof of $u\in C(\R^+,L^2(\R^3))$, where $u$ is a solution of $(NS_2)$ given by Friedrich approximation. We point out that we can use this method to show the same results in the case \cite{HK}.\\
$\bullet$ By inequality (\ref{eqth21}) we get
$$\limsup_{t\rightarrow0}\|u(t)\|_{L^2}\leq\|u^0\|_{L^2}.$$
Then, Proposition \ref{prop1}-(3) implies that $$\limsup_{t\rightarrow0}\|u(t)-u^0\|_{L^2}=0,$$
which ensures continuity at 0.\\
$\bullet$ Let $t_0>0$. For $\varepsilon\in(0,t_0/2)$ and $n\in\N$, put the following function
$$v_{n,\varepsilon}(t)=u_{\varphi(n)}(t+\varepsilon).$$
Applying the same method to prove the uniqueness to $u_{\varphi(n)}$ and $v_{n,\varepsilon}$, and using (\ref{eqf1}) we get
$$\|u_{\varphi(n)}(t+\varepsilon)-u_{\varphi(n)}(t)\|^{2}_{L^{2}}\leq\|u_{\varphi(n)}(\varepsilon)-u_{\varphi(n)}(0)\|^{2}_{L^{2}} \exp\Big(cF_n(t)\Big),$$
where
$$F_n(t)=\int_{0}^{t}(\|\partial_{3} \nabla_{h}u_{\varphi(n)}\|^{2}_{L^2}+\|\partial_{3}u_{\varphi(n)}\|^{2}_{L^2}+\|\nabla_{h}u_{\varphi(n)}\|^{2}_{L^2}).$$
By using inequalities (\ref{thapp1})-(\ref{thapp2}), we get
$$\begin{array}{lcl}F_n(t)&\leq&\displaystyle \|\partial_3u^0\|_{L^2}e^{a_{\alpha,\beta}t}+\|\partial_3u^0\|_{L^2}\Big(\int_{0}^{t}e^{a_{\alpha,\beta}z}dz\Big)+\frac{\|u^0\|^{2}_{L^2}}{2}\\
&\leq&\displaystyle 2\|\partial_3u^0\|_{L^2}e^{a_{\alpha,\beta}t}+\frac{\|u^0\|^{2}_{L^2}}{2}.
\end{array}$$
For $t\in[0,2t_0]$, we have
$$F_n(t)\leq 2\|\partial_3u^0\|_{L^2}e^{2a_{\alpha,\beta}t_0}+\frac{\|u^0\|^{2}_{L^2}}{2}=M_{\alpha,\beta}(t_0,u^0).$$
Then, for $t=t_0$ and $t=t_0-\varepsilon$, we get
\begin{equation}\label{eqf1}
\|u_{\varphi(n)}(t_0+\varepsilon)-u_{\varphi(n)}(t_0)\|^{2}_{L^{2}}\leq\|u_{\varphi(n)}(\varepsilon)-u_{\varphi(n)}(0)\|^{2}_{L^{2}} \exp\Big(cM_{\alpha,\beta}(t_0,u^0)\Big),
\end{equation}
\begin{equation}\label{eqf2}
\|u_{\varphi(n)}(t_0-\varepsilon)-u_{\varphi(n)}(t_0)\|^{2}_{L^{2}}\leq\|u_{\varphi(n)}(\varepsilon)-u_{\varphi(n)}(0)\|^{2}_{L^{2}} \exp\Big(cM_{\alpha,\beta}(t_0,u^0)\Big).
\end{equation}
The idea is to lower the terms on the left and increase the terms on the right of the inequalities (\ref{eqf1}) and (\ref{eqf2}).\\
For the right term, we write
 $$ \|u_{\varphi(n)}(\varepsilon)-u_{\varphi(n)}(0)\|^{2}_{L^{2}}=\|u_{\varphi(n)}(\varepsilon)\|^{2}_{L^{2}}+\|u_{\varphi(n)}(0)\|^{2}_{L^{2}}
 -2Re\langle u_{\varphi(n)}(\varepsilon),u_{\varphi(n)}(0)\rangle_{L^2}.$$
By using inequality (\ref{thapp1}) we obtain
 $$\begin{array}{lcl}\|u_{\varphi(n)}(\varepsilon)-u_{\varphi(n)}(0)\|^{2}_{L^{2}}
 &\leq&2\|u^0\|^{2}_{L^{2}}-2Re\langle u_{\varphi(n)}(\varepsilon),u_{\varphi(n)}(0)\rangle_{L^2}\\
 &\leq&2\|u^0\|^{2}_{L^{2}}-2Re\langle u_{\varphi(n)}(\varepsilon),u^0\rangle_{L^2}-2Re\langle u_{\varphi(n)}(\varepsilon),u_{\varphi(n)}(0)-u^0\rangle_{L^2}.
 \end{array}$$
But $$|\langle u_{\varphi(n)}(\varepsilon),u_{\varphi(n)}(0)-u^0\rangle_{L^2}|\leq\|u_{\varphi(n)}(\varepsilon)\|_{L^2}\|u_{\varphi(n)}(0)-u^0\|_{L^2}
\leq\|u^0\|_{L^2}\|u_{\varphi(n)}(0)-u^0\|_{L^2}$$
 then $$\lim_{n\rightarrow\infty}|\langle u_{\varphi(n)}(\varepsilon),u_{\varphi(n)}(0)-u^0\rangle_{L^2}|=0.$$
 On the other hand, and by using that $u_{\varphi(n)}(\varepsilon)$ converge weakly in $L^2(\R^3)$ to $u(\varepsilon)$, we get
 $$\liminf_{n\rightarrow\infty}\|u_{\varphi(n)}(\varepsilon)-u_{\varphi(n)}(0)\|^{2}_{L^{2}}\leq 2\|u^0\|^{2}_{L^{2}}-2Re\langle u(\varepsilon);u^0\rangle_{L^2}.$$
 For the left term, we have, for all $q,N\in\N$,
 $$\|J_N\Big(\theta_q.(u_{\varphi(n)}(t_0\pm\varepsilon)-u_{\varphi(n)}(t_0))\Big)\|^{2}_{L^2}
 \leq\|\theta_q.(u_{\varphi(n)}(t_0\pm\varepsilon)-u_{\varphi(n)}(t_0))\|^{2}_{L^2}\leq
 \|u_{\varphi(n)}(t_0\pm\varepsilon)-u_{\varphi(n)}(t_0)\|^{2}_{L^2}.$$
 Using (\ref{theq6}) we get, for $q$ big enough,
 $$\|J_N\Big(\theta_q.(u(t_0\pm\varepsilon)-u(t_0))\Big)\|^{2}_{L^2}
 \leq \liminf_{n\rightarrow\infty}\|u_{\varphi(n)}(t_0\pm\varepsilon)-u_{\varphi(n)}(t_0)\|^{2}_{L^2}.$$
 Then
$$\|J_N\Big(\theta_q.(u(t_0\pm\varepsilon)-u(t_0))\Big)\|^{2}_{L^2}
 \leq 2\Big(\|u^0\|^{2}_{L^{2}}-Re\langle u(\varepsilon);u^0\rangle_{L^2}\Big)\exp(cM_{\alpha,\beta}(t_0,u^0)).$$
By applying the Monotonic Convergence Theorem in the order $N\rightarrow\infty$ and $q\rightarrow\infty$, we get
$$\|u(t_0\pm\varepsilon)-u(t_0))\|^{2}_{L^2}
 \leq 2\Big(\|u^0\|^{2}_{L^{2}}-Re\langle u(\varepsilon);u^0\rangle_{L^2}\Big)\exp(cM_{\alpha,\beta}(t_0,u^0)).$$
Using the continuity at 0 and make $\varepsilon\rightarrow0$, we get the continuity at $t_0$, which ends the proof.

\end{document}